



\input amstex

\documentstyle{amsppt}
  \magnification=1100
  \hsize=6.2truein
  \vsize=9.0truein
  \hoffset 0.1truein
  \parindent=2em

\NoBlackBoxes

\newcount\mycitestyle \mycitestyle=1


\newcount\theTime
\newcount\theHour
\newcount\theMinute
\newcount\theMinuteTens
\newcount\theScratch
\theTime=\number\time
\theHour=\theTime
\divide\theHour by 60
\theScratch=\theHour
\multiply\theScratch by 60
\theMinute=\theTime
\advance\theMinute by -\theScratch
\theMinuteTens=\theMinute
\divide\theMinuteTens by 10
\theScratch=\theMinuteTens
\multiply\theScratch by 10
\advance\theMinute by -\theScratch

\def\today{{\number\day\space
 \ifcase\month\or
  January\or February\or March\or April\or May\or June\or
  July\or August\or September\or October\or November\or December\fi
 \space\number\year}}



















\define\biggnm#1{                            
  \bigg|\bigg|#1\bigg|\bigg|}

\define\bignm#1{                             
  \big|\big|#1\big|\big|}











\define\Cpx{\bold C}                         


\define\dif{\text{\it d}}                    







\define\eqdef{{\;\overset\text{def}\to=\;}}     







\define\fpamalg#1{{\dsize\;                  
     \operatornamewithlimits*_{#1}\;}}


\define\freeprod#1#2{\mathchoice             
     {\operatornamewithlimits{\ast}_{#1}^{#2}}
     {\raise.5ex\hbox{$\dsize\operatornamewithlimits{\ast}
      _{#1}^{#2}$}\,}
     {\text{oops!}}{\text{oops!}}}

\define\freeprodi{\mathchoice                
     {\operatornamewithlimits{\ast}
      _{\iota\in I}}
     {\raise.5ex\hbox{$\dsize\operatornamewithlimits{\ast}
      _{\sssize\iota\in I}$}\,}
     {\text{oops!}}{\text{oops!}}}

























\define\id{\text{\rm id}}                    




\define\Integers{\bold Z}                    




\define\Lambdao{{\Lambda\oup}}               




\define\lrnm#1{\left|\left|#1\right|\right|} 

\define\lspan{\text{\rm span}@,@,@,}         









\define\nm#1{||#1||}                         


\define\Naturals{{\bold N}}                  







\define\oup{^{\text{\rm o}}}                 

\define\owedge{{                             
     \operatorname{\raise.5ex\hbox{\text{$
     \ssize{\,\bigcirc\llap{$\ssize\wedge\,$}\,}$}}}}}

\define\owedgeo#1{{                          
     \underset{\raise.5ex\hbox
     {\text{$\ssize#1$}}}\to\owedge}}

\define\Pc{{\Cal P}}                         









\define\QED{\newline                         
            \line{$\hfill$\qed}\enddemo}









\define\Reals{{\bold R}}                     


\define\restrict{\lower .3ex                 
     \hbox{\text{$|$}}}




\define\smd#1#2{\underset{#2}\to{#1}}          

\define\smdb#1#2{\undersetbrace{#2}\to{#1}}    

\define\smdbp#1#2#3{\overset{#3}\to            
     {\smd{#1}{#2}}}

\define\smdbpb#1#2#3{\oversetbrace{#3}\to      
     {\smdb{#1}{#2}}}

\define\smdp#1#2#3{\overset{#3}\to             
     {\smd{#1}{#2}}}

\define\smdpb#1#2#3{\oversetbrace{#3}\to       
     {\smd{#1}{#2}}}

\define\smp#1#2{\overset{#2}\to                
     {#1}}                                     










\define\tocdots                              
  {\leaders\hbox to 1em{\hss.\hss}\hfill}    

\define\tr{\text{\rm tr}}                    




























  \newcount\bibno \bibno=0
  \def\newbib#1{\advance\bibno by 1 \edef#1{\number\bibno}}
  \ifnum\mycitestyle=1 \def\cite#1{{\rm[\bf #1\rm]}} \fi
  \def\scite#1#2{{\rm[\bf #1\rm, #2]}}


  \newcount\ignorsec \ignorsec=0
  \def\notasec{\ignorsec=1}

  \newcount\secno \secno=0
  \def\newsec#1{\procno=0 \subsecno=0 \ignorsec=0
    \advance\secno by 1 \edef#1{\number\secno}
    \edef\currentsec{\number\secno}}

  \newcount\subsecno
  \def\newsubsec#1{\procno=0 \advance\subsecno by 1 \edef#1{\number\subsecno}
    \edef\currentsec{\number\secno.\number\subsecno}}

  \newcount\appendixno \appendixno=0
  \def\newappendix#1{\procno=0 \ignorsec=0 \advance\appendixno by 1
    \ifnum\appendixno=1 \edef\appendixalpha{\hbox{A}}
      \else \ifnum\appendixno=2 \edef\appendixalpha{\hbox{B}} \fi
      \else \ifnum\appendixno=3 \edef\appendixalpha{\hbox{C}} \fi
      \else \ifnum\appendixno=4 \edef\appendixalpha{\hbox{D}} \fi
      \else \ifnum\appendixno=5 \edef\appendixalpha{\hbox{E}} \fi
      \else \ifnum\appendixno=6 \edef\appendixalpha{\hbox{F}} \fi
    \fi
    \edef#1{\appendixalpha}
    \edef\currentsec{\appendixalpha}}

  \newcount\procno \procno=0
  \def\newproc#1{\advance\procno by 1
   \ifnum\ignorsec=0 \edef#1{\currentsec.\number\procno}
   \else \edef#1{\number\procno}
   \fi}

  \newcount\tagno \tagno=0
  \def\newtag#1{\advance\tagno by 1 \edef#1{\number\tagno}}



\define\AoMn{{A\otimes M_n(\Cpx)}}

\define\BoMn{{B\otimes M_n(\Cpx)}}

\define\bpt#1{{                                
     \left[ #1 \right]_{\varphi\otimes\tr}}}   

\define\bpi#1{{                                
     \left[ #1 \right]_{\varphi\otimes\id}}}   

\define\iotn{{\id\otimes\tr_n}}

\define\ooMn{{1\otimes M_n(\Cpx)}}

\define\ooon{{1\otimes1_n}}

\define\poin{{\varphi\otimes\id_n}}

\define\potn{{\varphi\otimes\tr_n}}

\notasec
  \newtag{\AFP}
 \newproc{\GermainsThm}
  \newtag{\onto} 
  \newtag{\unitaries}
  \newtag{\unitarycond}
 \newproc{\OrderedKth}
 \newproc{\KnuclearInRep}
 \newproc{\InfProj}
  \newtag{\SigmaKK}
\newsec{\PrelimS}
 \newproc{\TwoFreeProj}
 \newproc{\uBuFree}
  \newtag{\varphiub}
 \newproc{\FreeEndom}
 \newproc{\GNSnotEnough}
\newsec{\FreenessMatrices}
 \newproc{\FreeMatrices}
  \newtag{\FreeFamily}
 \newproc{\EveryEntry}
  \newtag{\yEveryEntry}
  \newtag{\potneij}
\newsec{\EigenfreenessS}
 \newproc{\EigenfreenessD}
 \newproc{\SuffEigenfree}
 \newproc{\EigenfreeMn}
\newsec{\TracialS}
 \newproc{\KthTracial}
  \newtag{\KopTracialEF}
  \newtag{\scaleTracialEF}
 \newproc{\OrderedKthNonNucl}
\newsec{\NonTracialS}
 \newproc{\CompPosDef}
 \newproc{\CompPosProp}
 \newproc{\GetUnder}
 \newproc{\DphiNotInv}
 \newproc{\PropInf}
  \newtag{\Dauau}
  \newtag{\fkafkp}
 \newproc{\EqProj}
 \newproc{\EFPropInf}

\newbib{\AndersonBlackadarHaagerup}
\newbib{\Avitzour}
\newbib{\CuntzZZAddMult}
\newbib{\CuntzZZDimFunct}
\newbib{\CuntzZZKthCertain}
\newbib{\DykemaZZFreeProdR}
\newbib{\DykemaZZFaithful}
\newbib{\DykemaZZSimpl}
\newbib{\DykemaZZSNU}
\newbib{\DykemaHaagerupRordam}
\newbib{\DykemaRordamZZPI}
\newbib{\GermainZZKKeq}
\newbib{\GermainZZKthFFP}
\newbib{\GermainZZAmalgFIC}
\newbib{\RordamZZUHFII}
\newbib{\VilladsenZZPerf}
\newbib{\VoiculescuZZSymmetries}
\newbib{\VoiculescuZZMult}
\newbib{\VoiculescuZZCircSemiCirc}
\newbib{\VoiculescuZZRandMat}
\newbib{\VDNbook}

\topmatter

  \title Projections in free product C$^*$--algebras
  \endtitle

  \leftheadtext{Dykema, R\o rdam}
  \rightheadtext{Projections in free product C$^*$--algebras}

  \author Kenneth J\. Dykema and Mikael R\o rdam
  \endauthor

  \date \today \enddate

  \affil 
    Department of Mathematics and Computer Science \\
    Odense Universitet, Campusvej 55 \\
    DK-5230 Odense M \\
    Denmark
  \endaffil

  \address
    Department of Mathematics and Computer Science,
    Odense Universitet, Campusvej 55,
    DK-5230 Odense M,
    Denmark
  \endaddress

  \abstract
    Consider the reduced free product of C$^*$--algebras,
    $(A,\varphi)=(A_1,\varphi_1)*(A_2,\varphi_2)$, with respect to states
    $\varphi_1$ and
    $\varphi_2$ that are faithful.
    If $\varphi_1$ and $\varphi_2$ are traces, if the so--called Avitzour
    conditions are satisfied, (i.e\. $A_1$ and $A_2$ are not ``too small''
    in a specific sense) and if $A_1$ and $A_2$ are nuclear, then it is shown
    that the positive cone, $K_0(A)^+$, of the
    $K_0$--group of $A$ consists of those elements $g\in K_0(A)$
    for which $g=0$ or $K_0(\varphi)(g)>0$.
    Thus, the ordered group $K_0(A)$ is weakly unperforated.
        
    If, on the other hand, $\varphi_1$ or $\varphi_2$ is not a trace and if
    a certain condition weaker than the Avitzour conditions hold, then $A$
    is properly infinite.
  \endabstract


\endtopmatter

\document \TagsOnRight \baselineskip=18pt

\vskip3ex
\noindent{\bf Introduction and statement of the main results.}
\vskip3ex

The reduced free product of
C$^*$--algebras~\cite{\VoiculescuZZSymmetries}, \cite{\Avitzour}, (see
also the book~\cite{\VDNbook}), is the appropriate construction in
Voiculescu's theory of freeness~\cite{\VoiculescuZZSymmetries},
\cite{\VDNbook}.  Given unital C$^*$--algebras $A_1$ and $A_2$ with
states $\varphi_1$ and, respectively, $\varphi_2$, whose GNS
representations are faithful, we denote the corresponding reduced free
product by
$$ (A,\varphi)=(A_1,\varphi_1)*(A_2,\varphi_2). \tag{\AFP} $$ Recall
that $\varphi$ is a trace if and only if both $\varphi_1$ and
$\varphi_2$ are traces.  Moreover,~\cite{\DykemaZZFaithful}, $\varphi$
is faithful if and only if both $\varphi_1$ and $\varphi_2$ are
faithful.

The K-theory of $A$ can be calculated, at least when $A_1$ and $A_2$ are
nuclear, by the following theorem of Emmanuel Germain.

\proclaim{Theorem
\GermainsThm~(\cite{\GermainZZKKeq},~\cite{\GermainZZKthFFP})}
Let $A$ be a reduced free product C$^*$--algebra as in~(\AFP), and suppose that
$A_1$ and $A_2$ are nuclear.
Then there is an exact sequence of K--groups,
$$ \matrix
\Integers\cong&K_0(\Cpx)&\overset{(K_0(i_1),-K_0(i_2))}\to\longrightarrow
&K_0(A_1)\oplus K_0(A_2)&\overset{K_0(j_1)+K_0(j_2)}\to\longrightarrow
&K_0(A) \\ \vspace{2ex}
&\uparrow&&&&\downarrow \\ \vspace{2ex}
&K_1(A)&\overset{K_1(j_1)+K_1(j_2)}\to\longleftarrow
&K_1(A_1)\oplus K_1(A_2)&\overset{(K_1(i_1),-K_1(i_2))}\to\longleftarrow
&K_1(\Cpx)&=0,
\endmatrix $$
where
$i_k\colon\Cpx\to A_k$
is the unital $^*$--homomorphism and where
$j_k\colon A_k\to A$
is the unital embedding arising from the construction of the reduced free
product~(\AFP).
\endproclaim

It follows in particular that
$$ K_0(A) = K_0(j_1)(K_0(A_1))+K_0(j_2)(K_0(A_2)) \tag{\onto} $$
whenever $A_1$ and $A_2$ are nuclear. 

Let us consider the case where both $\varphi_1$ and $\varphi_2$
are traces. The pair $(A_1,\varphi_1)$, $(A_2,\varphi_2)$ is then said to
satisfy the Avitzour condition (c.f.\ \cite{\Avitzour}) if there exist
unitaries
$$ a\in A_1,\qquad b,c\in A_2 \tag{\unitaries} $$
which satisfy
$$ \varphi_1(a)=0,\qquad\varphi_2(b)=\varphi_2(c)=0=\varphi_2(b^*c).
\tag{\unitarycond} $$
Is this is the the case, then $A$ is simple and $\varphi$ is the unique tracial
state on $A$~(\cite{\Avitzour}).
It was proved in~\cite{\DykemaHaagerupRordam}
that under the same conditions, then $A$ is of stable rank one.

Let $P(A)$ denote the set of projections in $A$, set $P_n(A)=P(M_n(A))$,
and set $P_{\infty}(A)=\cup_{n=1}^{\infty}P_n(A)$.
Murray-von Neumann equivalence of projections yields an
equivalence relation on $P_{\infty}(A)$, denoted by $\sim$.
(If $p\in M_n(A)$ and $q \in M_m(A)$, then $p \sim q$ if there exists
$v\in M_{n,m}(A)$ such that $vv^* = p$ and $v^*v = q$.)
Moreover, if $p,q\in P_{\infty}(A)$, then $p \lesssim q$ if $p$ is equivalent
to a subprojection of $q$.

By the definition of the $K_0$-group of a (unital) C$^*$-algebra, we
have a map $[\,\cdot\,]_0\colon P_{\infty}(A)\to K_0(A)$, and $[p]_0=[q]_0$ if
$p \sim q$.
The image of this map is called {\it the positive cone}
of $K_0(A)$, and is denoted by $K_0(A)^+$.
The {\it scale} of $K_0(A)$ is defined to be the subset
$\{\,[p]_0\mid p\in P(A)\}$ of $K_0(A)^+$, and it is denoted by $\Sigma(A)$.
In the present
case, where $A$ is unital and has a faithful trace, $(K_0(A),K_0(A)^+)$
is an ordered abelian group (with $g\le h$ if $h-g\in K_0(A)^+$).
Notice that the trace $\varphi$ on $A$ induces a state
$K_0(\varphi) \colon K_0(A) \to {\Reals}$ given by
$K_0(\varphi)([p]_0-[q]_0)=\varphi(p)-\varphi(q)$.
Notice also that
$[p]_0 \le [q]_0$ if $p \lesssim q$.

The ordered group $K_0(A)$ is said to be {\it weakly unperforated} if
whenever $ng \ge 0$ for some positive integer $n$ and some $g\in K_0(A)$ it
follows that $g\ge 0$.

The order structure of the $K_0$-group of the Choi algebra and its
generalizations, ${C^*_r(\Integers_n*\Integers_m)}$, $n \ge 2$, $m \ge 3$,
was considered in \cite{\AndersonBlackadarHaagerup}.
They proved that
for each $t\in(\frac1n\Integers+\frac1m\Integers)\cap(0,1)$ there is a
projection in $C_r^*(\Integers_n*\Integers_m)$ of trace $t$.
They could not decide
whether all elements $g\in K_0(C_r^*(\Integers_n*\Integers_m))$ of positive
trace can be
realized by projections, (i.e\. are positive) or if
$K_0(C_r^*(\Integers_n*\Integers_m))$ is weakly unperforated.
Theorem~\OrderedKth{} below answers both of these
questions in the affirmative.

Jesper Villadsen \cite{\VilladsenZZPerf} has recently proven that not
all simple C$^*$-algebras have weakly unperforated $K_0$-groups.
This
shows that it will be hard (if not impossible) to obtain general
theorems about the positive cone of the $K_0$-group of arbitrary simple
C$^*$-algebras.
In the light of Villadsen's example and of our theorem
below, we conclude that the class of reduced free product C$^*$-algebras
(that we are considering) behaves rather nicely.

\proclaim{Theorem \OrderedKth}
Let $A$ be a reduced free product C$^*$--algebra as in~(\AFP), where
$\varphi_1$ and $\varphi_2$ are faithful traces.
Suppose that there are unitaries as
in~(\unitaries) which satisfy~(\unitarycond), and suppose that (\onto)
holds (which is the case if $A_1$ and $A_2$ are nuclear).
It follows
that 
$$ K_0(A)^+=\{g\in K_0(A)\mid K_0(\varphi)(g)>0\}\cup\{0\} $$
and
$$ \Sigma(A)=\{g\in K_0(A)\mid0<K_0(\varphi)(g)<1\}\cup\{0,1\}. $$
\endproclaim
It follows easily from this theorem that $K_0(A)$ is weakly
unperforated.
Moreover, 
$$K_0(\varphi)(K_0(A)) =
K_0(\varphi_1)(K_0(A_1))+K_0(\varphi_2)(K_0(A_2))\subseteq\Reals,$$
and if $G$ denotes the subgroup $\{g\in K_0(A)\mid K_0(\varphi)(g)=0\}$ of
$K_0(A)$, then
$$ K_0(A) = K_0(\varphi)(K_0(A)) \oplus G,$$
$$ K_0(A)^+ = \{(t,g) \in K_0(\varphi)(K_0(A)) \oplus G \mid t > 0\}
\cup \{(0,0)\}.$$ The group $G$ is, in other words, the
infinitesimal subgroup of $K_0(A)$.

It is a consequence of $A$ being of stable rank one, that the map
$[\,\cdot\,]_0 \colon P_{\infty}(A) \to K_0(A)$ induces an injection (and
hence a bijection)
$$ P_{\infty}(A)/\!\!\sim\quad\longrightarrow\;K_0(A)^+. $$
Under the conditions of Theorem~\OrderedKth{} this bijection is an order
isomorphism when $P_{\infty}(A)/\!\sim$ is equipped with the order
relation $\lesssim$ and when $K_0(A)^+$ is equipped with the algebraic
order relation (as above).
Thus, by using Theorem~\OrderedKth{} to find $K_0(A)^+$ and $\Sigma(A)$, we can
also classify Murray--von Neumann equivalence classes of projections in $A$ and
find their ordering.

\proclaim{Remark \KnuclearInRep}\rm
The assumption in Germain's Theorem above that
$A_1$ and $A_2$ be nuclear can be
loosened somewhat to the assumption that $A_k$ be K--nuclear
in the GNS representation of $\varphi_k$ $(k=1,2)$.
For the definition of this concept and the properties mentioned below,
see~\scite{\GermainZZAmalgFIC}{\S4}.
A C$^*$--algebra $A$ is K--nuclear in a given
representation $\pi$ if either
$A$ is nuclear or if $A=C^*_r(G)$ for a K--amenable group $G$ and $\pi$ is
the reduced regular representation.
Moreover, if $(A,\varphi)=(A_1,\varphi_1)*(A_2,\varphi_2)$ where each $A_j$ is
K--nuclear in the GNS representation of $\varphi_j$, then $A$ is
K--nuclear in the 
GNS representation of $\varphi$.
Thus, Theorem~\OrderedKth{} also applies to the case of the reduced
free product of finitely many C$^*$--algebras.
\endproclaim

Consider now the case where either $\varphi_1$ or $\varphi_2$ is not a
trace.
Then, as remarked above, $\varphi$ is not a trace.
If there are
unitaries as in (\unitaries) satisfying (\unitarycond) and if $a$ and
$b$ are in the centralizer of $\varphi_1$, respectively, $\varphi_2$,
then by Avitzour's theorem (\cite{\Avitzour}), $A$ is simple and does
not admit a tracial state. 

Recall that a projection $p$ in a C$^*$-algebra $A$ is said to be
infinite is it is equivalent to a proper subprojection of itself.
If
there exist two mutually orthogonal subprojections $p_1$ and $p_2$ of
$p$ such that $p \sim p_1 \sim p_2$, then $p$ is said to be {\it
properly infinite}.
A unital C$^*$-algebra $A$ is called infinite,
respectively, properly infinite, if its unit is infinite, respectively,
properly infinite.
Note that no infinite C$^*$--algebra has a faithful, tracial state and that no
properly infinite C$^*$--algebra has a tracial state.
A simple unital C$^*$-algebra is said to be {\it purely infinite}
if all its non-zero hereditary subalgebras contain an infinite
projection (see \cite{\CuntzZZKthCertain}).
It is an open problem whether every
simple unital
C$^*$-algebra, that does not admit a trace, must be purely infinite.
The
theorem below should be viewed in the light of that question.

\proclaim{Theorem \InfProj}
Let $A$ be a reduced free product C$^*$--algebra as in~(\AFP).
Suppose that $\varphi$ is faithful and is not a trace
and suppose there are unitaries as in (\unitaries) satisfying (\unitarycond).
Then $A$ is properly infinite.
Hence
$$ \Sigma(A)=K_0(A)^+=K_0(A). \tag{\SigmaKK} $$
\endproclaim
(Note that we assume less than Avitzour required, because we do not ask
unitaries to be in the centralizers of states.)

When proving these theorems, we define and make use of the intermediate concept
of {\it eigenfree} C$^*$--algebras.
Although the only C$^*$--algebras we prove are eigenfree are the reduced
free product C$^*$--algebras having unitaries as required in
Theorems~\OrderedKth{} and~\InfProj, we believe that more general reduced free
product C$^*$--algebras are eigenfree.
We have written our proofs in terms of eigenfreeness in order to allow for easy
use in possible generalizations.

\vskip3ex
\noindent{\bf\S\PrelimS. Preliminaries}
\vskip3ex

Basic to our investigation of projections in free product C$^*$--algebras is
the following result, which can be found in~\cite{\AndersonBlackadarHaagerup},
and can also be derived from Voiculescu's multiplicative free
convolution~\cite{\VoiculescuZZMult}.
See~\scite{\DykemaZZSimpl}{2.7} for a related result and
see~\cite{\DykemaZZSNU} for a detailed proof.

\proclaim{Proposition \TwoFreeProj}
Let $A$ be a C$^*$--algebra and $\varphi$ a faithful state on $A$ having
faithful GNS representation.
Suppose $p,q\in A$ are projections that are free with respect to $\varphi$.
If $\varphi(p)<\varphi(q)$ then $\nm{p(1-q)}<1$ and there is $v\in A$ such that
$v^*v=p$ and $vv^*\le q$.
\endproclaim

Another useful and well--known fact is the following.
(The proof is a particularly easy example of free etymology.)
As usual, we employ the notation, for subsets $S_1,S_2,\ldots,S_m$ of $A$,
$$ \Lambdao(S_1,S_2,\ldots,S_m)\eqdef\{a_1a_2a_3\cdots a_n
\mid n\in\Naturals,\,a_j\in S_{\iota_j},\,
\iota_1\ne\iota_2,\iota_2\ne\iota_3,\ldots,\iota_{n-1}\ne\iota_n\}. $$
\proclaim{Lemma \uBuFree}
Let $(A,\varphi)$ be a C$^*$--noncommutative probability space and let $B$ be a
subset of $A$.
Suppose $u\in A$ is a unitary,
that $\varphi(u)=0$ and that $B$ and $\{u,u^*\}$ are free.
Then $B$ and $u^*Bu$ are free.
\endproclaim
\demo{Proof}
We may suppose $B$ is a subalgebra of $A$ containing the unit of $A$.
By the freeness assumption, for $b\in B$, $\varphi(u^*bu)=0$ if and only if
$\varphi(b)=0$.
To show that $B$ and $u^*Bu$ are free it thus suffices to show that
$$ \varphi(b_1(u^*b_2u)b_3(u^*b_4u)\cdots b_{n-2}(u^*b_{n-1}u)b_n)=0
\tag{\varphiub} $$
whenever $n$
is odd, $n\ge3$, $b_j\in B$, $\varphi(b_j)=0$, $(2\le j\le n-1)$, and for $j=1$
and $j=n$ either $b_j=1$ or $\varphi(b_j)=0$.
However, the word $b_1u^*b_2ub_3\cdots u^*b_{n-1}ub_n$ is easily seen to equal
an alternating product in $\{b\in B\mid\varphi(b)=0\}$ and $\{u,u^*\}$.
By freeness of $B$ and $\{u,u^*\}$,~(\varphiub) holds.
\QED

\proclaim{Lemma \FreeEndom}
Let $I$ be an index set and let $(A_\iota,\varphi_\iota)$ be a
C$^*$--noncommutative probability space, $(\iota\in I)$, where each
$\varphi_\iota$ is faithful.
Let $(B,\psi)$ be a C$^*$--noncommutative probability space with $\psi$
faithful.
Let
$$ (A,\varphi)=\freeprodi(A_\iota,\varphi_\iota) $$
be the reduced free product C$^*$--algebra.
Given unital $^*$--homomorphisms, $\pi_\iota\colon A_\iota\to B$, such that
$\psi\circ\pi_\iota=\varphi_\iota$ and $(\pi_\iota(A_\iota))_{\iota\in I}$ is
free in $(B,\psi)$, there is a $^*$--homomorphism, $\pi\colon A\to B$ such that
$\pi\restrict_{A_\iota}=\pi_\iota$ and $\psi\circ\pi=\varphi$.
\endproclaim
\demo{Proof}
Let $B_0$ be the C$^*$--subalgebra of $B$ generated by
$\bigcup_{\iota\in I}\pi_\iota(A_\iota)$ and let $\psi_0=\psi\restrict_{B_0}$.
If the GNS representation associated to $\psi_0$ is faithful on $B_0$, then by
Voiculescu's construction~\cite{\VoiculescuZZSymmetries}, (see
also~\cite{\VDNbook}), $(B_0,\psi_0)$ is canonically
isomorphic to
$$ \freeprodi(A_\iota,\varphi_\iota) $$
via an isomorphism having the desired properties
The faithfulness of $\psi$, however, shows that $\psi_0$ is faithful, hence has
faithful GNS representation.
\QED

The following example shows that the hypothesis that $\psi$ be faithful is
essential, and cannot be replaced with the weaker hypothesis that the GNS
representation of $\psi$ be faithful.

\proclaim{Example \GNSnotEnough}\rm
Using the notation of~\cite{\DykemaZZSimpl}, let
$$ (B_1,\psi_1)\eqdef(\smdp\Cpx{3/4}p\oplus\smd\Cpx{1/4})*
(\smdp\Cpx{2/3}q\oplus\smd\Cpx{1/3}). $$
Let $B=M_2(B_1)$ and let $\psi$ be
the state on $M_2(B_1)$ given by
$\psi(\smallmatrix b_{11}&b_{12}\\b_{21}&b_{22}\endsmallmatrix)
=\psi_1(b_{11})$.
Although $\psi$ is not faithful, clearly the GNS representation of $\psi$ is
faithful on $B$.
Let
$$ \align
(A_1,\varphi_1)&=(\smdp\Cpx{3/4}r\oplus\smd\Cpx{1/4}), \\ \vspace{2ex}
(A_2,\varphi_2)&=(\smdp\Cpx{2/3}s\oplus\smd\Cpx{1/3}),
\endalign $$
and
$$ (A,\varphi)\eqdef(A_1,\varphi_1)*(A_2,\varphi_2). $$
Then $s$ is Murray--von Neumann equivalent in $A$ to a proper
subprojection of $r$.
Let $\rho_j\colon A_j\to B$ be the unital $^*$--homomorphisms such that
$\rho_1(r)=\left(\smallmatrix p&0\\0&0\endsmallmatrix\right)$ and
$\rho_2(s)=\left(\smallmatrix q&0\\0&1\endsmallmatrix\right)$.
Then $\psi\circ\rho_j=\varphi_j$.
But there cannot be a $^*$--homomorphism $\rho\colon A\to B$ such that
$\rho(r)=\rho_1(r)$ and $\rho(s)=\rho_2(s)$,
because $\rho_2(s)$ is not equivalent to a subprojection of
$\rho_1(r)$.
\endproclaim

\vskip3ex
\noindent{\bf\S\FreenessMatrices.  Matrices of free random variables.}
\vskip3ex

\proclaim{Theorem \FreeMatrices}
Let $n\in\Naturals$ and suppose $(A,\varphi)$ is a $^*$--noncommutative
probability
space having random variables $x_{ij}\in A$, $(1\le i\le j\le n)$ and a unital
subalgebra $B\subseteq A$ such that
\roster
\item"(i)" $x_{ii}$ is a semicircular element with
$\varphi(x_{ii}^2)=1$, $(1\le i\le n)$;
\item"(ii)" $x_{ij}$ is a circular element with
$\varphi(x_{ij}^*x_{ij})=1$, $(1\le i<j\le n)$;
\item"(iii)" the family of sets of random variables,
$$ \bigl(B,(\{x_{ii}\})_{1\le i\le n},
(\{x_{ij}^*,x_{ij}\})_{1\le i<j\le n}\bigr) \tag{\FreeFamily} $$
is free.
\endroster
Let $(e_{ij})_{1\le i,j\le n}$ be a system of matrix units for $M_n(\Cpx)$.
Consider the noncommutative probability space
$(A\otimes M_n(\Cpx),\varphi\otimes\tr_n)$, where $\tr_n$ is the tracial state
on $M_n(\Cpx)$.
Consider also the random variable
$$ x=\frac1{\sqrt n}\biggl(\sum_{1\le i\le n}x_{ii}\otimes e_{ii}
+\sum_{1\le i<j\le n}(x_{ij}\otimes e_{ij}+x_{ij}^*\otimes e_{ji})\biggr). $$
Then $x$ is semicircular and $\{x\}$ and $\BoMn$ are free.
\endproclaim
\demo{Proof}
Voiculescu, using his matrix model~\cite{\VoiculescuZZRandMat},
proved in~\cite{\VoiculescuZZCircSemiCirc}
that $x$ is semicircular.
Moreover, from~\cite{\DykemaZZFreeProdR} it follows that $\{x\}$ and
$1\otimes M_n(\Cpx)$ are free.

Let $1_n$ denote the unit of $M_n(\Cpx)$ and let $\id_n$ denote the identity
map on $M_n(\Cpx)$.
For $p\in\Naturals$ and $a\in\AoMn$, let
$$ \align
\bpt{x^p}&=x^p-(\potn)(x^p)\cdot(\ooon), \\
\bpi a&=a-(\poin)(a).
\endalign $$

\proclaim{Lemma \EveryEntry}
Let $m\in \Naturals$ and $p_1,p_2,\ldots,p_m\in\Naturals$.
Let $d_0,d_1,\ldots,d_m\in\ooMn$ and assume that
$(\iotn)(d_j)=0$ whenever $1\le j\le m-1$.
Let
$$ y=d_0\bpt{x^{p_1}}d_1\bpt{x^{p_2}}\cdots d_{m-1}\bpt{x^{p_m}}d_m.
\tag{\yEveryEntry} $$
Then
$(\poin)(y)=0$.
\endproclaim
Note that $(\potn)(y)=0$ by freeness of $\{x\}$ and $\ooMn$.
The lemma gives more, namely that every matrix entry of $y$ has zero
expectation.
\demo{Proof}
To show that $(\poin)(y)=0$ it will suffice to show that
$$ (\potn)\bigl((1\otimes e_{1i})y(1\otimes e_{j1})\bigr)=0. \tag{\potneij} $$
whenever $1\le i,j\le n$.
But absorbing $1\otimes e_{1i}$ into $d_0$ and $1\otimes e_{j1}$ into $d_m$,
we see that $(1\otimes e_{1i})y(1\otimes e_{j1})$ is a word having the same
form as $y$.
As mentioned above, the freeness of $\{x\}$ and $\ooMn$ then
implies~(\potneij).
\QED

Continuing with the proof of Theorem~\FreeMatrices, take arbitrary
$m\in\Naturals$, $p_1,p_2,\ldots,p_m\in\Naturals$ and
$a_0,a_1,\ldots,a_{m}\in\BoMn$ such that
$(\potn)(a_j)=0$ for every $1\le j\le m-1$.
To prove the theorem it will suffice to show that $(\potn)(z)=0$, where
$$ z=a_0\bpt{x^{p_1}}a_1\bpt{x^{p_2}}\cdots
a_{m-1}\bpt{x^{p_m}}a_{m}. $$
Writing $a_j=\bpi{a_j}+(\poin)(a_j)$ and distributing,
we write $z$ as a sum of $2^{m+1}$ terms,
$$ \align
z&=\bpi{a_0}\bpt{x^{p_1}}\bpi{a_1}\bpt{x^{p_2}}\cdots
\bpi{a_{m-1}}\bpt{x^{p_m}}\bpi{a_{m}}+ \\ \vspace{1ex}
&+\bpi{a_0}\bpt{x^{p_1}}\bpi{a_1}\bpt{x^{p_2}}\cdots
\bpi{a_{m-1}}\bpt{x^{p_m}}(\poin)(a_{m})+ \\ \vspace{1ex}
&+\qquad\cdots\qquad+ \\
&+(\poin)(a_0)\bpt{x^{p_1}}(\poin)(a_1)\bpt{x^{p_2}}\cdots
(\poin)(a_{m-1})\bpt{x^{p_m}}\bpi{a_{m}}+ \\ \vspace{1ex}
&+(\poin)(a_0)\bpt{x^{p_1}}(\poin)(a_1)\bpt{x^{p_2}}\cdots \\
&\qquad\qquad\qquad\qquad\qquad\qquad\qquad\qquad\qquad
\cdots(\poin)(a_{m-1})\bpt{x^{p_m}}(\poin)(a_{m}).
\endalign $$
Since each $(\poin)(a_j)\in\ooMn$ and since if $1\le j\le m-1$ we
have $(\potn)\circ(\poin)(a_j)=0$, it follows from the freeness of $\{x\}$ and
$\ooMn$ that $\potn$ of the last term is zero.
Each of the remaining $2^{m+1}-1$ terms is of the form
$$ t=b_0\bpi{a_{j(1)}}b_1\bpi{a_{j(2)}}\cdots b_{k-1}\bpi{a_{j(k)}}b_k $$
where $1\le k\le m$, $0\le j(1)<j(2)<\cdots<j(k)\le m$ and
$$ b_0=\cases\ooon&\text{if }j(1)=0\\
(\poin)(a_0)\bpt{x^{p_1}}\cdots(\poin)(a_{j(1)-1})\bpt{x^{p_{j(1)}}}
&\text{if }j(1)>0,\endcases $$
if $1\le\,l\le k-1$ then
$$ b_l=\bpt{x^{p_{j(l)+1}}}(\poin)(a_{j(l)+1})\bpt{x^{p_{j(l)+2}}}\cdots
(\poin)(a_{j(l+1)-1})\bpt{x^{p_{j(l+1)}}} $$
and
$$ b_k=\cases\bpt{x^{p_{j(k)+1}}}(\poin)(a_{j(k)+1})\bpt{x^{p_{j(k)+2}}}\cdots
\\ \qquad\qquad\qquad\cdots
(\poin)(a_{m-1})\bpt{x^{p_{m}}}(\poin)(a_{m})&\text{if }j(k)<m
\\ \vspace{1ex}
\ooon&\text{if }j(k)=m.\endcases $$
From Lemma~\EveryEntry, we see that $(\poin)(b_l)=0$ whenever $0\le l\le k$,
except when $l=0$ and $j(0)=0$ or when $l=k$ and $j(k)=m$, (i.e\. when
$b_0=\ooon$ or $b_k=\ooon$).
Thus, excepting the cases just mentioned, every matrix entry of each $b_l$
belongs to the $*$-subalgebra of $A$ generated by
$\{x_{ij}\mid1\le i\le j\le n\}$ and evaluates to $0$ under $\varphi$.
In addition, every matrix entry of each $\bpi{a_j}$ belongs to $B$ and
evaluates to zero under $\varphi$.
Hence it follows from the freeness of~(\FreeFamily) that every matrix entry of
$t$ evaluates to zero under $\varphi$.
Therefore, each $(\potn)(t)=0$ and hence $(\potn)(z)=0$.
\QED

\vskip3ex
\noindent{\bf\S\EigenfreenessS. Eigenfreeness}
\vskip3ex

\proclaim{Definition \EigenfreenessD}\rm
Let $A$ be a unital C$^*$--algebra with state $\varphi$.
We say that $(A,\varphi)$ is {\it eigenfree} if there is a unital
$^*$--endomorphism, $\rho$, of $A$ and a Haar unitary with respect to
$\varphi$, $u\in A$, such that
$\rho(A)$ and $\{u\}$ are $^*$--free in $(A,\varphi)$ and
$\varphi\circ\rho=\varphi$.
\endproclaim

\proclaim{Proposition \SuffEigenfree}
Let $A_1\ne\Cpx$ and $A_2\ne\Cpx$ be unital C$^*$--algebras having
faithful states $\varphi_1$ and, respectively, $\varphi_2$.
Let
$$ (A,\varphi)=(A_1,\varphi_1)*(A_2,\varphi_2) $$
be the C$^*$--algebra reduced free product.
If there are unitaries
$a\in A_1$ and $b,c\in A_2$ such that $\varphi_1(a)=0$ and
$\varphi_2(b)=\varphi_2(c)=0=\varphi_2(b^*c)$
then $(A,\varphi)$ is eigenfree by an endomorphism $\rho$ such that
$K_0(\rho)\colon K_0(A)\to K_0(A)$ restricts to the identity map on the image
of $K_0(A_j)\to K_0(A)$, $(j=1,2)$.
\endproclaim
\demo{Proof}
Let $\rho_k$, $(k=1,2)$, be the $^*$--homomorphisms
$$ \rho_k\colon A_k\to A $$
defined by
$\rho_1(x)=a^*bxb^*a$ and $\rho_2(y)=ba^*yab^*$.
Then each $\rho_k$ is injective.
Moreover, since $\varphi_1(x)=0$ implies $\varphi(a^*bxb^*a)=0$ and
$\varphi_2(y)=0$ implies $\varphi(ba^*yab^*)=0$, we have
$\varphi\circ\rho_k=\varphi_k$, $(k=1,2)$.
Let us show that $\rho_1(A_1)$ and $\rho_2(A_2)$ are free.
It will suffice to show that $\varphi(z)=0$ whenever
$z\in\Lambdao(\rho_1(A_1\oup),\rho_2(A_2\oup))$.
But for such a word, $z$, no cancellation occurs and we see that $z$ is equal
to an element of $\Lambdao(A_1\oup,A_2\oup)$, so by freeness, $\varphi(z)=0$.
By~\cite{\DykemaZZFaithful}, the state $\varphi$ is faithful.
It then follows from Lemma~\FreeEndom{} that there is a $^*$--endomorphism,
$\rho$, of $A$ such that
$\rho\circ j_k=\rho_k$, $(k=1,2)$ and $\varphi\circ\rho=\varphi$.

Let $u=a^*cac^*$.
Clearly $u$ is a Haar unitary in $(A,\varphi)$.
We now show that $\{u\}$ and $\rho(A)$ are $^*$--free in $(A,\varphi)$, which
will complete the proof that $(A,\varphi)$ is eigenfree.
By the freeness of $\rho_1(A_1)$ and $\rho_2(A_2)$,
$\lspan\Lambdao(\rho_1(A_1\oup),\rho(A_2\oup))$ is dense in $\rho(A)\oup$.
Therefore, it will suffice to show that $\varphi(z)=0$ whenever
$$ z\in\Lambdao\bigl(a^*b(A_1\oup)b^*a,\,ba^*(A_2\oup)ab^*,\,
\{u^n\mid n\in\Naturals\}\cup\{(u^*)^n\mid n\in\Naturals\}\bigr). $$
Expand each $u$ as $a^*cac^*$ and each $u^*$ as $ca^*c^*a$.
Now the only cancellations which may occur are
$$ \align
(a^*bxb^*a)(a^*cac^*)&=a^*bx(b^*c)ac^* \\
(ca^*c^*a)(a^*bxb^*a)&=ca^*(c^*b)xb^*a \\
(a^*cac^*)(ba^*yab^*)&=a^*ca(c^*b)a^*yab^* \\
(ba^*yab^*)(ca^*c^*a)&=ba^*ya(b^*c)a^*c^*a,
\endalign $$
for $x\in A_1\oup$ and $y\in A_2\oup$.
Making these cancellations, $z$ is seen to be equal to an element of
$\Lambdao(A_1\oup,A_2\oup)$, so by freeness $\varphi(z)=0$.

Since $\rho$ on the copies of $A_j$ in $A$ is conjugation by a unitary, we
easily see that $K_0(\rho)$ is the identity map on the image of
$K_0(A_j)\to K_0(A)$, $(j=1,2)$.
\QED

\proclaim{Proposition \EigenfreeMn}
Let $A$ be a C$^*$--algebra with state $\varphi$ and let $n\in\Naturals$.
If $(A,\varphi)$ is eigenfree by an endomorphism $\rho$ and some
unitary then
$(A\otimes M_n(\Cpx),\varphi\otimes\tr_n)$ is eigenfree by the
endomorphism
$\rho^{n^2}\otimes\id_n$, where $\tr_n$ is the tracial state on
$M_n(\Cpx)$ and
$\id_n$ is the identity map on $M_n(\Cpx)$.
\endproclaim
\demo{Proof}
Let $\rho$ and $u$ be as in Definition~\EigenfreenessD.
Because $\varphi\circ\rho=\varphi$ and using~\scite{\VDNbook}{2.5.5(iii)}, the
$n^2$ unitaries,
$$ u,\,\rho(u),\,\rho^2(u),\,\ldots,\,\rho^{n^2-1}(u), $$
are Haar unitaries and the family
$$ (\rho^{n^2}(A),\,\{u\},\,\{\rho(u)\},\,\ldots,\,\{\rho^{n^2-1}(u)\}) $$
is $^*$--free.
Using the continuous functional calculus, we find semicircular elements
$x_{ii}\in A$ such that $\varphi(x_{ii}^2)=1$, $(1\le i\le n)$, and circular
elements $x_{ij}\in A$ such that $\varphi(x_{ij}^*x_{ij})=1$,
$(1\le i<j\le n)$, such that
$$ \bigl(\rho^{n^2}(A),
(\{x_{ii}\})_{1\le i\le n},(\{x_{ij}\})_{1\le i<j\le n}\bigr) $$
is $^*$--free.
Let $(e_{ij})_{1\le i,j\le n}$ be a system of matrix units for $M_n(\Cpx)$ and
let
$$ x=\frac1{\sqrt n}\biggl(\sum_{1\le i\le n}x_{ii}\otimes e_{ii}
+\sum_{1\le i<j\le n}(x_{ij}\otimes e_{ij}+x_{ij}^*\otimes e_{ji})\biggr). $$
Then by Theorem~\FreeMatrices, $x$ is semicircular and
$\rho^{n^2}(A)\otimes M_n(\Cpx)$ and $\{x\}$ are free.
Using the continuous functional calculus we obtain a Haar unitary,
$v\in C^*(\{x\})$, which then satisfies that $\rho^{n^2}(A)\otimes M_n(\Cpx)$
and $\{v\}$ are $^*$--free.

Thus the endomorphism $\rho^{n^2}\otimes\id_n$ of $A\otimes M_n(\Cpx)$ and the
Haar unitary $v$ give 
that $(A\otimes M_n(\Cpx),\potn)$ is eigenfree.
\QED

\vskip3ex
\noindent{\bf\S\TracialS.  The tracial case}
\vskip3ex

In this section we prove Theorem~\OrderedKth,
which follows from the proposition below, in conjunction with
Proposition~\SuffEigenfree.

\proclaim{Proposition \KthTracial}
Let $A$ be a unital C$^*$--algebra with faithful, tracial state $\tau$ and
suppose that $(A,\tau)$ is eigenfree by an endomorphism $\rho$.
Suppose that $G$ is a subgroup of $K_0(A)$ on which $K_0(\rho)$ is the identity
map.
Then
$$ G\cap K_0(A)^+=\{x\in G\mid K_0(\tau)(x)>0\}\cup\{0\} \tag{\KopTracialEF} $$
and
$$ G\cap\Sigma(A)=\{x\in G\mid0<K_0(\tau)(x)<1\}\cup\{0,1\}
\tag{\scaleTracialEF} $$
\endproclaim
\demo{Proof}
Since $\tau$ is a faithful trace on $A$, the inclusion $\subseteq$ is clear
in both~(\KopTracialEF) and~(\scaleTracialEF).
To show $\supseteq$ in~(\KopTracialEF), let $x\in G$ have $K_0(\tau)(x)>0$.
Then there is $n\in\Naturals$ and there are projections
$p,q\in\AoMn$ such
that $x=[p]-[q]$.
We now use Proposition~\EigenfreeMn{}; let $\rho_n$ denote the endomorphism
$\rho^{n^2}\otimes\id_n$ and let $u_n\in\AoMn$ denote the Haar unitary such
that $\rho_n(\AoMn)$ and $\{u_n\}$ are $^*$--free.
Then
$x=[\rho_n(p)]-[\rho_n(q)]$.
Thus $(\tau\otimes\tr_n)(\rho_n(p))>(\tau\otimes\tr_n)(\rho_n(q))$.
Moreover, by Lemma~\uBuFree, $\rho_n(p)$ and $u_n^*\rho_n(q)u_n$
are free and clearly the traces of $\rho_n(q)$ and  $u_n^*\rho_n(q)u_n$ are the
same.
So by Proposition~\TwoFreeProj, $u_n^*\rho_n(q)u_n$, and thus also $\rho_n(q)$,
is equivalent in $\AoMn$ to a subprojection, say $r$, of $\rho_n(p)$.
Then $x=[\rho_n(p)-r]\in K_0(A)^{+}$.
This proves~(\KopTracialEF).

The truth of~(\scaleTracialEF) now follows similarly.
Indeed, if $x\in G$ satisfies $0<K_0(\tau)(x)<1$ then by~(\KopTracialEF)
there is $n\in\Naturals$ and a projection $p\in\AoMn$ such that
$x=[p]=[\rho_n(p)]$.
Let $e=1\otimes e_{11}\in\AoMn$, where $e_{11}$ is a rank--one projection in
$M_n(\Cpx)$.
Then since $K_0(\tau)(x)<1$ we have
$(\tau\otimes\tr_n)(p)<(\tau\otimes\tr_n)(e)$.
Clearly $e=\rho_n(e)$.
By the same argument as was applied above in the proof of~(\KopTracialEF),
we see
that $\rho_n(p)$ is equivalent to a subprojection, say $s$, of $e$.
Thus $x=[s]\in\Sigma(A)$.
\QED

In fact, the above proposition together with Proposition~\SuffEigenfree{}
proves the following result, which is more general than Theorem~\OrderedKth.
\proclaim{Proposition \OrderedKthNonNucl}
Let $A$ be a reduced free product C$^*$--algebra as in~(\AFP), where
$\varphi_1$ and $\varphi_2$ are faithful traces.
Suppose that there are unitaries as
in~(\unitaries) which satisfy~(\unitarycond).
Consider the subgroup
$$ G=K_0(j_1)(K_0(A_1))+K_0(j_2)(K_0(A_2))\subseteq K_0(A). $$
Then
$$ G\cap K_0(A)^+=\{x\in G\mid K_0(\varphi)(x)>0\}\cup\{0\} $$
and
$$ G\cap\Sigma(A)=\{x\in G\mid0<K_0(\varphi)(x)<1\}\cup\{0,1\}. $$
\endproclaim

\vskip3ex
\noindent{\bf\S\NonTracialS.  The non--tracial case.}
\vskip3ex

In this section we will prove Theorem~\InfProj.

We will make use of the comparison theory for positive elements in a
C$^*$--algebra that was introduced by
J\. Cuntz~\cite{\CuntzZZAddMult},~\cite{\CuntzZZDimFunct} (see
also~\cite{\RordamZZUHFII}), and which we describe in the definition and
proposition below.

\proclaim{Definition \CompPosDef}\rm
Let $A$ be a unital C$^*$--algebra with positive cone denoted $A^+$ and let
$a,b\in A^+$.
Write
$$ a\lesssim b $$
if there is a sequence, $(x_n)_{n=1}^\infty$, in $A$ such
that $\lim x_n^*bx_n=a$.

If $\varphi$ is a state on $A$, define the function,
$D_\varphi\colon A^+\to[0,1]$ by
$$ D_\varphi(a)=\lim_{\epsilon\searrow0}\varphi(f_\epsilon(a)), $$
where
$$ f_\epsilon(t)=\cases 0&\text{if } 0\le t\le\epsilon \\
(t-\epsilon)/\epsilon&\text{if }\epsilon\le t\le2\epsilon \\
1&\text{if }t\ge2\epsilon. \endcases $$
(If $\varphi$ were a trace then $D_\varphi$ would be a dimension function.)
\endproclaim

\proclaim{Proposition \CompPosProp{} (\cite{\CuntzZZAddMult},
\cite{\CuntzZZDimFunct}, \cite{\RordamZZUHFII})}
Let $A$ be a unital C$^*$--algebra.
Then
\roster
\item"(i)" the relation $\lesssim$ on $A$ is transitive and reflexive;
\item"(ii)" $\lesssim$, when restricted to the projections of $A$, gives the
usual Murray--von~Neumann ordering.
\endroster
Let $a,b\in A^+$ and $x\in A$.
Then
\roster
\item"(iii)" if $f_\epsilon(a)\lesssim b$ for every $\epsilon>0$, then
$a\lesssim b$;
\item"(iv)" if $a\le b$ then $a\lesssim b$;
\item"(v)" if $f\colon \Reals_+\to\Reals_+$ is a continuous function with
$f(0)=0$ then $f(a)\lesssim a$;
\item"(vi)" $xx^*\lesssim x^*x$;
\item"(vii)" $D_\varphi(p)=\varphi(p)$ if $p\in A$ is a projection;
\item"(viii)" $\varphi(a)=\int_0^{\nm a}D_\varphi(f_t(a))\dif t$.
\endroster
\endproclaim

\proclaim{Lemma \GetUnder}
Let $A$ be a unital C$^*$--algebra and let $\varphi$ be a faithful state on
$A$.
Suppose that $p$ is a projection in $A$, $a\in A^+$ and $p$ and $a$ are free
with respect to $\varphi$.
\roster
\item"(i)" If $D_\varphi(a)<\varphi(p)$ then $a\lesssim p$.
\item"(ii)" If $\varphi(p)<D_\varphi(a)$ then $p\lesssim a$.
\endroster
\endproclaim
\demo{Proof}
Let $A''$ denote the von Neumann algebra generated by the image of $A$ under
the GNS representation of $\varphi$, and denote also by $\varphi$ the
normal extension of $\varphi$ to $A''$.

For~(i), it will suffice to show that $f_\epsilon(a)\lesssim p$ for every
$\epsilon>0$.
Set
$$ q=\chi_{[\epsilon,\infty)}(a)\in A''. $$
Then $f_\epsilon(a)\le q\le f_{\epsilon/2}(a)$, $f_\epsilon(a)q=f_\epsilon(a)$,
and $p$ and $q$ are free with respect to $\varphi$.
Therefore
$$ \varphi(q)\le\varphi(f_{\epsilon/2}(a))\le D_\varphi(a)<\varphi(p), $$
and hence, by Proposition~\TwoFreeProj, $\nm{q(1-p)}<1$.
Setting $\lambda=1-\nm{q(1-p)}^2$, it follows that
$$ q(1-p)q\le\nm{q(1-p)}^2q=(1-\lambda)q, $$
so  $qpq\ge\lambda q$.
Hence
$$ \lambda f_\epsilon(a)=\lambda f_\epsilon(a)^{1/2}qf_\epsilon(a)^{1/2}
\le f_\epsilon(a)^{1/2}qpqf_\epsilon(a)^{1/2}
=f_\epsilon(a)^{1/2}pf_\epsilon(a)^{1/2}, $$
which implies that
$$ f_\epsilon(a)\lesssim \lambda f_\epsilon(a)
\le f_\epsilon(a)^{1/2}pf_\epsilon(a)^{1/2}\lesssim pf_\epsilon(a)p\le p. $$

For~(ii), choose $\epsilon>0$ such that $\varphi(f_\epsilon(a))>\varphi(p)$.
Set
$$ q=\chi_{[\epsilon,\nm a]}(a)\in A''. $$
Then, as above, $p$ and $q$ are free with respect to $\varphi$ and
$f_\epsilon(a)\le q\le f_{\epsilon/2}(a)$, so
$\varphi(p)<\varphi(q)$, whence $\nm{p(1-q)}<1$.
Setting $\lambda=1-\nm{p(1-q)}^2$ we get
$\lambda p\le pqp\le pf_{\epsilon/2}(a)p$ which gives that
$$ p\lesssim\lambda p\le pf_{\epsilon/2}(a)p
\lesssim f_{\epsilon/2}^{1/2}(a)pf_{\epsilon/2}^{1/2}(a)
\le f_{\epsilon/2}(a)\lesssim a. $$
\QED

\proclaim{Lemma \DphiNotInv}
Let $A$ be a unital C$^*$--algebra and let $\varphi$ be a state on $A$ that is
not a trace.
Then there are $a\in A^+$ and a unitary $u\in A$ such that
$$ D_\varphi(a)<D_\varphi(u^*au). $$
\endproclaim
\demo{Proof}
Suppose for contradiction that $D_\varphi(a)=D_\varphi(u^*au)$ for every
$a\in A^+$ and unitary $u\in A$.
Then by Proposition~\CompPosProp(vii),
$$ \varphi(a)=\int_0^{\nm a}D_\varphi(f_t(a))\dif t
=\int_0^{\nm a}D_\varphi(u^*f_t(a)u)\dif t
=\int_0^{\nm a}D_\varphi(f_t(u^*au))\dif t=\varphi(u^*au) $$
for every $a\in A^+$ and every unitary $u\in A$.
But then $\varphi$ would be a trace.
\QED

\proclaim{Lemma \PropInf}
Let $A$ be a unital C$^*$--algebra and $\varphi$ a faithful state on $A$.
Suppose that $(A,\varphi)$ is eigenfree by the endomorphism $\rho$ and unitary
$v$, and suppose that $\varphi$ is not a trace.
Then for every $n$ large enough, the identity is a properly infinite projection
in $A\otimes M_n(\Cpx)$.
\endproclaim
\demo{Proof}
By Lemma~\DphiNotInv{} there are $a\in A^+$ and a unitary $u\in A$
such that $D_\varphi(a)<D_\varphi(u^*au)$.
Let $(e_{ij})_{1\le i,j\le n}$ be a system of matrix
units for $M_n(\Cpx)$ and for $1\le l\le n$ let
$f_l=\sum_{j=1}^le_{jj}$.
Let $k,n\in\Naturals$ be such that
$$ D_\varphi(a)<\frac kn<\frac{k+1}n<D_\varphi(u^*au). $$
Then in $\AoMn$ and for every $m\in\Naturals$
$$ \aligned
D_\potn(\rho^m(a)\otimes1_n)=&D_\varphi(a)<(\potn)(1\otimes f_k)< \\
&<(\potn)(1\otimes f_{k+1})<D_\varphi(u^*au)
=D_\potn(\rho^m(u^*au)\otimes1_n). 
\endaligned \tag{\Dauau} $$
By Proposition~\EigenfreeMn, $(\AoMn,\potn)$ is eigenfree by the endomorphism
$\rho^{n^2}\otimes\id_n$ and a unitary $w\in\AoMn$.
Thus, by Lemma~\uBuFree, $w^*(\rho^{n^2}(A)\otimes M_n(\Cpx))w$ and
$\rho^{n^2}(A)\otimes M_n(\Cpx)$ are free.
Hence, from~(\Dauau) and Lemma~\GetUnder,
$w^*(1\otimes f_{k+1})w\lesssim\rho^{n^2}(u^*au)\otimes1_n$
and
$\rho^{n^2}(a)\otimes1_n\lesssim w^*(1\otimes f_k)w$, which shows
that in $\AoMn$
$$ 1\otimes f_{k+1}\lesssim \rho^{n^2}(a)\otimes1_n
\lesssim1\otimes f_k. \tag{\fkafkp} $$
From this one easily obtains that $1\otimes f_n\lesssim1\otimes f_l$ in $\AoMn$
whenever $k\le m,l\le n$.
This again shows that $1\otimes f_k$ is properly infinite in $\AoMn$
and, furthermore, the identity of $A\otimes M_m(\Cpx)$ is properly infinite
whenever $m\ge k$.
\QED

\proclaim{Lemma \EqProj}
Let $A$ be a unital C$^*$--algebra and $\varphi$ a faithful state on $A$.
Suppose that $(A,\varphi)$ is eigenfree by an endomorphism $\rho$ and some
unitary.
Suppose that $\varphi$ is not a trace.
Let $n\in\Naturals$ be such that $\ooon\in\AoMn$ is properly infinite.
If $p,q\in\rho^{n^2}(A)\otimes M_n(\Cpx)$ are nonzero projections then
$p\lesssim q$ and $q\lesssim p$ in $\AoMn$.
\endproclaim
\demo{Proof}
Since $\ooon$ is properly infinite, there is a sequence, $(r_j)_{j=1}^\infty$,
of mutually orthogonal projections in $\AoMn$, each equivalent to $\ooon$.
Let $\rho_n$ denote the $^*$--endomorphism $\rho^{n^2}\otimes\id_n$ of $\AoMn$
and let $\varphi_n$ denote the state $\potn$ of $\AoMn$.
Since $\rho_n(\ooon)=\ooon$, we find a
sequence $(p_j)_{j=1}^\infty$ of projections in
$\rho_n(\AoMn)$ such that each $p_j$ is equivalent in
$\rho_n(\AoMn)$ to $p$ and $p_j\le\rho_n(r_j)$.
Since $\lim_{j\to\infty}\varphi_n(p_j)=0$, by replacing $p$ by some $p_j$ if
necessary, we may assume without loss of generality that
$\varphi_n(p)<\varphi_n(q)$.
By Lemma~\EigenfreeMn, there is a Haar unitary $w_n\in\AoMn$ such that
the image of $\rho_n$ and $\{w_n\}$ are $^*$--free.
By Lemma~\uBuFree, $p$ and $w^*_nqw_n$ are free.
Moreover, $\varphi_n(w^*_nqw_n)=\varphi_n(q)$.
By Proposition~\TwoFreeProj, $p\lesssim w_n^*qw_n$, so $p\lesssim q$.
By symmetry we get also $q\lesssim p$.
\QED

\proclaim{Lemma \EFPropInf}
Let $A$ be a unital C$^*$--algebra and $\varphi$ a faithful state on $A$.
Suppose that $(A,\varphi)$ is eigenfree and that $\varphi$ is not a trace.
Then the unit, $1$, of $A$ is properly infinite.
\endproclaim
\demo{Proof}
Using Lemma~\PropInf, let $n\in\Naturals$ be such that $\ooon$ is properly
infinite.
Let $\rho$ be a $^*$--endomorphism by which $(A,\varphi)$ is eigenfree and let
$\rho_n$ denote the $^*$--endomorphism $\rho^{n^2}\otimes\id_n$ of $\AoMn$.
Identify the unit, $1\in A$, with $e=1\otimes e_{11}\in\AoMn$,
where $e_{11}\in M_n(\Cpx)$ is a minimal projection.
By Lemma~\EqProj, $\ooon$ is equivalent to a subprojection, $r$, of $e$, which
is thus properly infinite.
Since $r$ is properly infinite, it has orthogonal subprojections,
$r_1,r_2,r_3,r_4$ such that each $r_j$ is equivalent to $r$.
Replacing $r$ with $\rho_n(r)$ and each $r_j$ with $\rho_n(r_j)$ if necessary,
we may assume without loss of generality that
$r,r_1,r_2,r_3,r_4\in\rho_n(\AoMn)$.
Then, by Lemma~\EqProj, $e-r\lesssim r_3$ and $e-r\lesssim r_4$.
Hence $e=r+(e-r)\lesssim r_1+r_3$ and $e\lesssim r_2+r_4$.
This shows that $e$, and thus the unit, $1$, of $A$, is properly infinite.
\QED

Having proved that $A$ is properly infinite,~(\SigmaKK) now follows
from~\scite{\CuntzZZKthCertain}{1.4}, since the set, $\Pc$, of all properly
infinite projections in $A$ is seen to satisfy the conditions ($\Pi_1$),
($\Pi_2$), ($\Pi_3$) and ($\Pi_4$) of~\cite{\CuntzZZKthCertain}.
Therefore, Theorem~\InfProj{} follows from Proposition~\SuffEigenfree{} and
Lemma~\EFPropInf.

\enddocument